\documentclass{amsart}
\usepackage{epsf,amssymb,amsfonts}
\usepackage{epsfig}
\usepackage[mathscr]{eucal}
\def\tauv{{\large $\tau$}$\!_{\mathrm{\scriptscriptstyle V}}$}
\def\taue{{\large $\tau$}$\!_{\mathrm{\scriptscriptstyle E}}$}
\def\tauee{{\large $\tau$}$\!_{\mathrm{\scriptscriptstyle E_{2}}}$}
\def\taup{{\large $\tau$}$\!_{\mathrm{\scriptscriptstyle P}}$}
\def\taum{{\large $\tau$}$\!_{\mathrm{\scriptscriptstyle M}}$}

\newtheorem{thm}{Theorem}
\newtheorem{lem}{Lemma}
\newtheorem{cor}{Corollary}
\newtheorem{exmp}{Example}
\newtheorem{defn}{Definition}

\numberwithin{equation}{section}
\DeclareMathOperator{\diam}{diam}
\DeclareMathOperator{\dm}{d}
\DeclareMathOperator{\p}{Pr}

\begin{document}
\bibliographystyle{plain}
\title[End compactifications in non-locally-finite graphs]{End compactifications in non-locally-finite graphs}

\author[B.\ Kr\"on]{B.\ Kr\"on$^{\textstyle{\,\star}}$}
\begin{abstract}
There are different definitions of ends in non-locally-finite graphs which
are all equivalent in the locally finite case. We prove the compactness of the
end-topology that is based on the principle of removing finite sets of vertices and give a proof of the
compactness of the end-topology that is constructed by the 
principle of removing finite sets of edges. For the latter case there exists
already a proof in \cite{cartwright93martin}, which only works on graphs with countably infinite vertex sets and in
contrast to which we do not use the Theorem of Tychonoff. We also construct a new topology of ends that
arises from the principle of removing sets of vertices with finite diameter and give applications that
underline the advantages of this new definition.
\end{abstract}

\thanks{$^{\textstyle{\star}}$ The author is supported by the START-project Y96-MAT of the
  Austrian Science Fund. Current address: Institut f\"ur Mathematik C, Technische Universit\"at Graz, Steyrergasse 30, 8010 - Graz, tel.: +43/316/873-4509, e-mail: kroen@finanz.math.tu-graz.ac.at}

\maketitle

\section{Introduction}
Ends of graphs can be seen as the directions along which sequences of vertices can tend to infinity.
Freudenthal \cite{freudenthal31enden} was the first who considered ends on a class of topological spaces, which we nowadays
call locally finite graphs. Usually ends are defined as equivalence classes of rays. Although one can find
various literature on ends of locally 
finite graphs, for an introduction to this topic see \cite{moeller95groups} or \cite{moeller92ends1}, there does not even exist a
standard definition of ends in the non-locally-finite case. Halin 
\cite{halin64unendliche} calls two rays equivalent, when both have infinitely many vertices in
common with a third ray. We call the arising equivalence classes \emph{vertex-ends}, on which the topological
considerations of Polat in \cite{polat79aspects1}, \cite{polat79aspects2} and \cite{polat90topological} are based. Cartwright, Soardi, Woess \cite{cartwright93martin}, Dicks and Dunwoody \cite{dicks89groups},
\cite{dunwoody82cutting} and Stallings \cite{stallings71group} prefer ends that are constructed by the
principle of removing finite sets of edges. We call them \emph{edge-ends}. For the latter approach to the
subject Dunwoody has proved the existence of so-called structure trees for graphs with more than one end,
which can be powerful tools in describing structures of infinite graphs with a strong action of its
automorphism group \cite{dunwoody82cutting}. An improved proof can be found in \cite{dicks89groups}. At the other hand each
vertex-end is contained in an edge-end, which means that vertex-ends can describe structures that cannot be
seen by using edge-ends. Cf. Example 3, Example 5 and Lemma 7. But none of these two approaches could yield convincing
arguments for taking one of them as the standard.  Our new construction of so-called \emph{metric ends}, that is only based on the principle of distinguishing between sets of finite and sets of infinite
diameters, enables us to describe structures that cannot even be seen by using vertex-ends (cf.\ Example 4, Example 5 and Section 8), but we do not claim to have reached the end of this discussion.\par
We want to see ends as points in the boundaries of the compactification of topologies on the set of vertices of a graph, first of all to obtain a concept for the convergence of sequences of vertices.\par
Throughout this article we will only use rather simple theorems of General Topology. For the most part we will use pure graphtheoretic arguments. We could see end compactifications in a more general topological context, but in our experience this is in most of the cases not too useful.\par
The \emph{edge-compactification} with corresponding topology \taue, which is the compactification that arises from the definition of edge-ends, is the only one from which we can construct a Hausdorff topology \tauee\ in a natural way, cf. \cite{cartwright93martin}, further explanations thereto can be found at the end of Section 4. Although it is a rather weak topology, it is the concept that is used by the majority of the authors.\par
Into this topology we can embed the topology \tauv\ of the \emph{vertex-compactification}, the counterpart of the vertex-ends.\par
The \emph{proper metric topology} \taup, arising from the definition of metric ends, is the strongest of the three topologies and has nearly all `good' topological properties but when we want it to be compact we have to add two artificial extra points (improper ends) and the resulting \emph{metric topology} \taum\ is not even $T_{0}$ although every sequence has a convergent subsequence with a unique limit. The reason for that is that there is no natural concept of local convergence in compactifications of discrete metric spaces which allows all open balls with positive real radius to be open. We can embed the other end compactifications into the proper metric topology but not into the metric compactification with its two extra points.\par
There do not arise any problems from these two extra points, when we want to use this topology for example to describe convergence of a transient random walk to a boundary (cf. Section 8), because the probability of the event that they occur as a limit is usually zero. The proper metric topology has furthermore some pretty topological properties concerning the study of 
quasi-isometries (cf.\ Theorem 6 and Theorem 7).\par
Random walks on graphs can tend to infinity along directions that cannot be described in a pure graph theoretical way without using the properties of the transition densities. It would be a subject of further researches to find homeomorphisms between graphtheoretical end compactifications and probability theoretical compactifications like the Martin compactification under corresponding preliminary restrictions to the random walk. In \cite{cartwright93martin} the Martin boundary is embedded into the space of edge-ends.\par
Figure 2 shows a synopsis of some topological properties of the end topologies that will be discussed in this article.\par
\newpage
{\scriptsize
\begin{tabular*}{7cm}{ccccc}

{\normalsize \taum} & & & & {\normalsize \tauee}\vspace{1mm}\\
compact, & & & & comp., Lind.,\\
Lindel\"of, $\neg T_{0}$ & & & & normal \\
\vspace{17mm}&\hspace*{1.5mm}&&\hspace*{15mm}&\\
{\normalsize \taup}&cont.\ surj.&{\normalsize \tauv}&cont.\ surj.&{\normalsize \taue}\vspace{1mm}\\

$\neg$comp., $\neg$Lind., &&comp., Lind., $T_{0}$,&&comp., Lind.,\\
Tychonoff,&&Def.\ of Halin && $\neg T_{0}$\\
if loc.\ count.\ fin.&&&&\\
then normal&&&&\\

\end{tabular*}}
\par
\vskip3mm
\setlength{\unitlength}{0.8mm}
\begin{picture}(1,1)
\put(15,29){\vector(0,1){14}}
\put(115,29){\vector(0,1){14}}
\put(31.5,19){\vector(1,0){15}}
\put(85,19){\vector(1,0){15}}
\put(18,35){{\scriptsize + 2 improper ends}}
\put(118,35){{\scriptsize + improper vertices}}
\put(135.5,30.78){{\scriptsize cf.\ \cite{cartwright93martin}}}
\end{picture}
\begin{center}
{\em Figure 2}
\end{center}
\section{A simple property of non-locally-finite graphs}
To give a slight impression of the properties of non-locally-finite graphs we first of all want to state a simple lemma which characterises graphs with infinite diameter. Throughout this article let $X=(VX,EX)$ be a connected graph without loops and multiple edges. A set of
vertices $e$ is called connected, if any two vertices in $e$ can be connected
by a path in $X$ that does not leave $e$. When we consider topological connectedness we will mention it explicitly if it is not clear from the context. We write
$e^{*}$ for the complement $VX\backslash e$ of $e$ and $\diam_{X}$ for the diameter with respect to the natural graph metric $\dm_{X}$ of $X$. The set of all connected components of $e^{*}$ is denoted by $\mathcal{C}(e)$. For the set of those connected components in $e^{*}$ that have a finite diameter we write $\mathcal{C}_{0}(e)$. A \emph{star ball} is a ball $K$ for which
$$\sup\{\diam C\mid C \in \mathcal{C}_{0}(K)\}=\infty.$$
A \emph{ray} is a sequence $(x_{n})_{n\in \mathbb{N}}$ of pairwise disjoint
vertices such that $x_{n}\sim x_{n+1}$ for all $n$.
\begin{lem}
\mbox{}\vspace{-4pt}
\begin{enumerate}
\item  A locally finite graph has infinite diameter if and only if it contains a ray.
\item The diameter of a non-locally-finite graph is infinite if and only if it contains a ray with infinite diameter or a star ball.
\end{enumerate}
\end{lem}

\begin{exmp}
Let\/ $\{P_{n}\mid n\in \mathbb{N}, \diam P_{n}=n \}$ be a set of disjoint paths. By joining together the initial vertices of these paths we obtain a graph\/ $X$ with a vertex\/ $x$ of infinite degree. The diameter of $X$ is infinite. There does not exist a ray with infinite diameter but for any natural radius\/ $r$ the ball\/ $B(x,r)$ is a star ball. See Figure 1.
\end{exmp}

\begin{center}

\setlength{\unitlength}{1mm}
\begin{picture}(160,60)
\put(50,20){\circle*{1}}
\put(41,32){\circle*{1}}
\put(50,35){\circle*{1}}
\put(50,50){\circle*{1}}
\multiput(59,32)(9,12){3}{\circle*{1}}
\multiput(65,25)(15,5){4}{\circle*{1}}
\multiput(50,20)(12,-2){6}{\line(6,-1){6}}

\put(35,31){$P_{1}$}
\put(43,50){$P_{2}$}
\put(69,55){$P_{3}$}
\put(107,43){$P_{4}$}
\put(48,16){$x$}
\put(50,20){\line(-3,4){9}}
\put(50,20){\line(0,1){30}}
\put(50,20){\line(3,4){27}}
\put(50,20){\line(3,1){60}}

\end{picture}
\vskip-0.4cm
{\em Figure 1}
\end{center}

\begin{proof}[Proof of Lemma 1]
The first statement is well known. To prove the second part of the lemma we have to show that there exists a ray with infinite diameter in any graph $X$ that does not contain a star ball.\par
If there is no star ball in $X$ then the complement of any vertex $x$ in $X$ must contain a connected component $e_{0}$ with infinite diameter. At least one of the components in the ball $K(x,1)^{*}$ which are contained in $e_{0}$ must have infinite diameter. Otherwise $K(x,1)$ would be a star ball. We call this component $e_{1}$. By induction we obtain a strictly decreasing sequence $(e_{r})_{r\in \mathbb{N}}$ of components of $K(x,r)^{*}$ with infinite diameter. In the proof of Lemma 3 we will see that there must exist a ray that has infinitely many vertices in common with every set $e_{r}$. Since the ray is not contained in any ball it must have infinite diameter.
\end{proof}

\section{The construction of vertex- and edge-topology}

We define
the \emph{vertex-boundary} $\theta e$ as the set of vertices in $e^{*}$ which
are adjacent to a vertex in $e$. $I\theta e:= \theta e^{*}$ is called
\emph{inner vertex-boundary} of $e$. The \emph{edge-boundary} $\delta e$ is defined as the
set of edges connecting vertices in $e$ with vertices in $e^{*}$. We call a
nonempty set of vertices $e$
\emph{vertex-} or \emph{edge-cut} if $\theta e$ or $\delta e$ are finite, respectively.\par
Note that every \emph{edge-cut} is also a \emph{vertex-cut}. The reversal of
this statement is not true in the general case, but in locally finite graphs
\emph{vertex-cut} and \emph{edge-cut} are equivalent terms. It will turn out that this is the reason why
the two topologies we are interested in are identical in the locally finite case.\par
A ray \emph{lies} in a set of vertices $e$ or is \emph{contained} in $e$, if $e$ contains all but finitely many elements of the ray. Sometimes we will use the terms \emph{contain} and \emph{lie} at the same time in the above sense as well as in the sense of set theoretic inclusion. A set of vertices $e$ \emph{separates} two rays, if one of them lies in $e$ and the other lies in $e^{*}$.\par
To rays are called \emph{vertex-equivalent} or \emph{edge-equivalent} if they cannot be separated by \emph{vertex-} or \emph{edge-cuts}, respectively. It is easy to see that these relations are
equivalence relations. Their equivalence classes are called \emph{vertex-}
and \emph{edge-ends} of $X$, respectively. As every edge-cut is a vertex-cut,
all vertex-ends are subsets of edge-ends. An end \emph{lies} in a set of vertices $e$ or is \emph{contained} in $e$, if all of its rays lie in $e$. The set of vertex-ends that lie in
$e$ is denoted by $\Omega_{\mathrm{\scriptscriptstyle V}}e$, the set of edge-ends in $e$
by $\Omega_{\mathrm{\scriptscriptstyle E}}e$. We write $\Omega_{\mathrm{\scriptscriptstyle
V/E}}X$ instead of $\Omega_{\mathrm{\scriptscriptstyle V/E}}VX$. Indeed, a vertex-end or
edge-end $\omega$ lies in a vertex-cut or edge-cut $e$ if and only if one of its rays lies in $e$,
respectively. $\omega$ lies either in $e$ or in $e^{*}$. We say that a
finite set of vertices $e$ \emph{separates} two ends, if they lie in different connected components of
$e^{*}$.

\begin{exmp}
The two sided infinite `ladder graph' has two ends. In locally finite graphs we do not have to distinguish between different types of ends. 
\end{exmp}

\begin{center}

\setlength{\unitlength}{1mm}
\begin{picture}(95,20)
\multiput(15,0)(15,0){5}{\line(0,1){15}}
\put(15,0){\vector(-1,0){15}}
\put(15,0){\vector(1,0){75}}
\put(15,15){\vector(-1,0){15}}
\put(15,15){\vector(1,0){75}}
\multiput(15,0)(15,0){5}{\circle*{1}}
\multiput(15,15)(15,0){5}{\circle*{1}}
\end{picture}\par\vskip0.5mm
{\em Figure 3}
\end{center}
\vskip0.5mm
\begin{lem}
For a graph\/ $X$ the set\/ $B_{\mathrm{\scriptscriptstyle V}}X:=\{e\cup
\Omega_{\mathrm{\scriptscriptstyle V}}e \mid e \subset VX \mbox{\ and\ \,} |\theta e| <
\infty \} $ is closed under finite intersection.
\end{lem}
\begin{proof}
For two sets $e_{1}\cup \Omega_{\mathrm{\scriptscriptstyle V}}e_{1}$ and
$e_{2}\cup \Omega_{\mathrm{\scriptscriptstyle V}}e_{2}$ in $B_{\mathrm{\scriptscriptstyle V}}X$ we have

$$(e_{1}\cup \Omega_{\mathrm{\scriptscriptstyle V}}e_{1}) \cap (e_{2}\cup
\Omega_{\mathrm{\scriptscriptstyle V}}e_{2})=(e_{1}\cap e_{2})\cup(\Omega_{\mathrm{\scriptscriptstyle
  V}}e_{1}\cap \Omega_{\mathrm{\scriptscriptstyle V}}e_{2}).$$
$e_{1}\cap e_{2}$ is a vertex-cut and $\Omega_{\mathrm{\scriptscriptstyle
  V}}e_{1}\cap \Omega_{\mathrm{\scriptscriptstyle V}}e_{2}$ the set of ends whose
rays are completely contained in $e_{1}$ and $e_{2}$ from some index onwards. Thus
$\Omega_{\mathrm{\scriptscriptstyle V}}e_{1}\cap \Omega_{\mathrm{\scriptscriptstyle
  V}}e_{2}=\Omega_{\mathrm{\scriptscriptstyle V}}(e_{1}\cap e_{2})$, and finally
$$(e_{1}\cup \Omega_{\mathrm{\scriptscriptstyle V}}e_{1}) \cap (e_{2}\cup
\Omega_{\mathrm{\scriptscriptstyle V}}e_{2})=(e_{1}\cap
e_{2})\cup\Omega_{\mathrm{\scriptscriptstyle V}}(e_{1}\cap e_{2}).$$
\end{proof}

We define the set $B_{\mathrm{\scriptscriptstyle E}}X$ analogously, for which the proof
of Lemma 2 can be copied word by word. Hence $B_{\mathrm{\scriptscriptstyle V}}X$ and
$B_{\mathrm{\scriptscriptstyle E}}X$ are bases of topological spaces $(VX\cup
\Omega_{\mathrm{\scriptscriptstyle V}}X,$\tauv$X)$ and $(VX\cup
\Omega_{\mathrm{\scriptscriptstyle E}}X,$\taue$ X)$, whose topologies \tauv$X$ and
\taue$X$ are called {\em vertex-topology\/} and {\em edge-topology}, respectively.

\section{Compactness of the vertex- and edge-topologies}

Let $K(z,n)$ denote the ball $\{x \in VX \mid \dm_{X}(x,z)\le n \}$ where
$\dm_{X}$ is the natural geodesic metric of the graph $X$.

\begin{lem}
If there is a vertex\/ $z$ for every strictly decreasing sequence\/ $(e_{n})_{n\in \mathbb{N}}$ of
connected edge- or vertex-cuts, such that\/ $K(z,n)$ is
a subset of\/ $e_{n}^{*}$, then there exists a ray\/ $L$, which has infinitely many
vertices in common with all cuts\/ $e_{n}$. In other words: The end of\/ $L$ lies
in all cuts\/ $e_{n}$.
\end{lem}

\begin{proof}
As $e_{1}$ is connected and a superior set of $e_{2}$ any vertex $x_{1}$ in $I\theta e_{1}$ can be connected by a path $\pi_{1}$ of vertices in
$e_{1}\backslash e_{2}$ with a vertex $x_{2}$ in $I\theta e_{2}$. By induction
we get a path $\pi_{n}$ such that the initial vertex of $\pi_{n}$ is adjacent to the last vertex of $\pi_{n-1}$ and $\pi_{n}$ is contained in $e_{n}\backslash e_{n+1}$ for every natural $n$ greater than one. The union $L$ of these paths must have finite intersection with every ball with center $z$ (eq. every bounded subset of $VX$). Thus $L$ is infinite and therefore it constitutes a ray. The corresponding end lies in all cuts $e_{n}$.
\end{proof}

\begin{lem}
Let\/ $e$ be a vertex-cut and\/ $\xi$ an infinite sequence of pairwise different
elements of\/ $e\cup \Omega_{\mathrm{\scriptscriptstyle V}}e$. If there is no
connected component of\/ $e$ containing infinitely many elements of the
sequence then\/ $\xi$ has an accumulation point in\/ $\theta e$.
\end{lem}

\begin{proof}
Under the given preliminaries there exists a vertex $x$ in $\theta e$ with
infinite degree, which is adjacent to infinitely many connected
components of $e$ containing elements of $\xi$. Every element $f\cup
\Omega_{\mathrm{\scriptscriptstyle V}}f$ of the base
$B_{\mathrm{\scriptscriptstyle V}}X$ which contains $x$ must contain
almost all of these components, since the vertex-boundary of $e$ is finite. The
same holds for every neighbourhood of $x$. Thus $x$ is an accumulation point of $\xi$.
\end{proof}

\begin{lem}
$(VX\cup \Omega_{\mathrm{\scriptscriptstyle V}}X,$\tauv$ X)$ is sequentially
compact. 
\end{lem}

\begin{proof}
Let $\xi$ be a sequence of pairwise different elements of $VX\cup
\Omega_{\mathrm{\scriptscriptstyle V}}X$. If $VX$ is the only vertex-cut, then
every element of $VX\cup \Omega_{\mathrm{\scriptscriptstyle E}}X$ is an
accumulation point of the sequence. Otherwise let $z$ be a vertex in the
complement of a vertex-cut $e$, in which lie infinitely many elements of
$\xi$. If the sequence has no accumulation point in $\theta e$, then by Lemma 4
there exists a connected component $e_{0}$ of $e$, such that infinitely
many elements of $\xi$ lie in $e_{0}$.\par
In the case that $\xi$ has no accumulation point in $\theta e_{0}$, we construct a connected vertex-cut
$e_{1}$, that is a subset of $e_{0}\backslash I\theta e_{0}$ and contains infinitely many elements of
$\xi$:\par
For a vertex $y$ in $\theta e_{0}$ with infinite degree, which is no accumulation point of $\xi$, there
exists a vertex-cut $f(y)$, that contains $y$ and for which $f(y)\cup \Omega_{\mathrm{\scriptscriptstyle
V}}f(y)$ contains only finitely many elements of $\xi$. Almost all vertices in $I\theta e_{0}$, that are
neighbours of $y$, lie in $f(y)$. Thus $e_{0}\backslash f(y)$ is a cut, such that infinitely many
elements of the sequence $\xi$ lie in it. This and the finiteness of $\theta e_{0}$ imply that
$$e_{1}^{\mathrm{\scriptscriptstyle A}}:=e_{0}\backslash \bigcup \{f(y)\mid y \in \theta e_{0},
\deg(y)=\infty\}$$
is a vertex-cut, which is a subset of $e_{0}$ and contains almost all elements of $I\theta e_{0}$. Thus
$$e_{1}^{\mathrm{\scriptscriptstyle B}}:=e_{1}^{\mathrm{\scriptscriptstyle A}}\backslash I\theta e_{0}$$
is also a vertex-cut. It is a subset of $e_{0}\backslash I\theta e_{0}$. If $\xi$ has no accumulation point
in $\theta e_{1}^{\mathrm{\scriptscriptstyle B}}$, once more by Lemma 3, there exists a connected component
$e_{1}$ of $e_{1}^{\mathrm{\scriptscriptstyle B}}$, that has again the properties requested before.\par
By induction we obtain a strictly decreasing sequence $(e_{n})_{n\in \mathbb{N}}$ of connected vertex-cuts,
in all of which lie infinitely many elements of the sequence $\xi$. Since $e_{n+1}$ is a subset of
$e_{n}\backslash I\theta e_{n}$ and the vertex $z$ lies in the complement of $e_{0}$ we have
$$K(z,n)\cap e_{n}=\emptyset $$
and thus
$$\bigcap_{n\in \mathbb{N}} e_{n}=\emptyset.$$
Following Lemma 3 we obtain a ray $L$, whose end $\omega$ lies in every cut $e_{n}$. Every neighbourhood $U$
of $\omega$ contains a base element $f\cup \Omega_{\mathrm{\scriptscriptstyle V}}f$ with $\omega \in
\Omega_{\mathrm{\scriptscriptstyle V}}f$, for which $\theta f$ is a subset of one of the balls $K(z,n)$. This
implies, that almost all cuts $e_{n}$ and thus infinitely many elements of the sequence $\xi$ lie in any
neighbourhood of $\omega$.
\end{proof}

\begin{lem}
$(VX\cup \Omega_{\mathrm{\scriptscriptstyle V}}X,$\tauv$X)$ is a Lindel\"of space.
\end{lem}

\begin{proof}
We may assume that an open cover $\mathcal{U}$ consists of base elements. Given a set $e\cup \Omega_{\mathrm{\scriptscriptstyle V}}e$ of the cover, we
choose a vertex $z$ in $e$. Since $\theta e$ is finite there exists a finite subcover
$\mathcal{U}_{1}$
of $e \cup \theta e$ in $\mathcal{U}$. The set $e_{1}$ of vertices in $\bigcup\mathcal{U}_{1}$ is a
vertex-cut. Every ray in $e_{1}$ must lie in one of those vertex-cuts $f$ for which $f\cup
\Omega_{\mathrm{\scriptscriptstyle V}}f$ is an element of the cover $\mathcal{U}_{1}$. The
same holds for any end containing this ray and therefore
$$e_{1}\cup \Omega_{\mathrm{\scriptscriptstyle V}}e_{1}=\bigcup\mathcal{U}_{1}.$$
By induction we obtain a sequence of vertex-cuts $(e_{n})_{n\in \mathbb N}$ and a corresponding finite
subcover $\mathcal{U}_{n}$ of $\mathcal{U}$, such that $K(z,n)$ is a subset of $e_{n}$. In other
words:
$$\qquad\ e_{n}\cup \Omega_{\mathrm{\scriptscriptstyle V}}e_{n}=\bigcup\mathcal{U}_{n} \mbox{\quad and \quad}\bigcup\{e_{n}\mid n\in \mathbb N\}=VX.$$
For an end ${\omega}$ and a vertex $z$ we define  $\dm_{\mathcal{\, U}}(\omega,z)$, the {\em distance of
$\omega$ to $z$ with respect to the cover $\mathcal{U}$}, as the minimal radius $r$ for which 
there exists a base element $f\cup \Omega_{\mathrm{\scriptscriptstyle V}}f$ in $\mathcal{U}$ such that
$\theta f$ is a subset of $K(z,r)$. Note that this radius exists for all vertex-ends. Let $f\cup
\Omega_{\mathrm{\scriptscriptstyle V}}f$ be a base element in $\mathcal{U}$. If the vertex 
boundary of $f$ is contained in $e_{n}$, then for an $x$ in $\theta e_{n}\cap f$ all connected components
of $(e_{n}\cup\theta e_{n})^{*}$ which are adjacent to $x$ are completely contained in $f$. Let
$\mathcal{W}_{n}$ denote the set of elements of $\mathcal{U}$ whose vertex-boundaries are subsets of $e_{n}$. We now choose a finite subcover $\mathcal{V}_{n}$ of
$\mathcal{W}_{n}$ which covers $\theta e_{n}\cap \bigcup \mathcal{W}_{n}$. Now $\mathcal{V}_{n}$ covers all
ends $\omega$ with $\dm_{\mathcal{\, U}}(\omega,z)\le n$ and the union
$$\bigcup_{n\in \mathbb N}\mathcal{U}_{n}\cup\mathcal{V}_{n}$$
is a countable covering of both $VX$ and $\Omega_{\mathrm{\scriptscriptstyle V}}X$ which meets the statement
of the lemma. 
\end{proof}

\begin{thm}
$(VX\cup \Omega_{\mathrm{\scriptscriptstyle V}}X,$\tauv$X)$ is a compact\/ $T_{0}$-space.
\end{thm}
 
\begin{proof}
A sequentially compact Lindel\"of space is compact. Thus we can conclude from Lemma 5 and Lemma 6, that $(VX\cup\Omega_{\mathrm{\scriptscriptstyle V}}X,$\tauv$X)$ is compact.\par
It remains to show that at least one of two elements $x$ and $y$ of $VX\cup
\Omega_{\mathrm{\scriptscriptstyle V}}X$ has an open environment, in which the other element is not
contained. We distinguish:

\begin{enumerate}
\item One of the two elements is a vertex. If $x$ is a vertex, then $(VX\backslash \{x\} \cup
\Omega_{\mathrm{\scriptscriptstyle V}}X)$ is the requested neighbourhood of $y$.
\item $x$ and $y$ are ends. In the complement of a finite set of vertices that separates $x$ and
$y$ there exist two disjoint base elements, of which one contains the end $x$ and the other $y$.
\end{enumerate}
\end{proof}
\begin{lem}
There exists a continuous surjection\/ $f$ from\/ $(VX\cup \Omega_{\mathrm{\scriptscriptstyle V}}X,$\tauv$X)$
onto\/ $(VX\cup \Omega_{\mathrm{\scriptscriptstyle E}}X,$\taue$X)$, whose restriction on\/ $VX$ is the
identity.
\end{lem}

\begin{proof}
We define the $f$-image of an end in $\Omega_{\mathrm{\scriptscriptstyle V}}X$ as the end in
$\Omega_{\mathrm{\scriptscriptstyle E}}X$ that contains it. The preimage of every base element of the
edge-topology is open in the vertex-topology.
\end{proof}

\begin{thm}
$(VX\cup \Omega_{\mathrm{\scriptscriptstyle E}}X,$\taue$X)$ is compact.
\end{thm}

\begin{proof}
The statement is a consequence of Theorem 1 and Lemma 7.
\end{proof}

\begin{exmp}
We connect two adjacent vertices\/ $x_{1}$ and\/ $x_{2}$ with each vertex of two disjoint rays\/ $L_{1}$ and\/ $L_{2}$, respectively. The resulting graph\/ $X_{1}$ has two vertex- and two edge-ends, see Figure 4a. When we take away the edge\/ $\{x_{1},x_{2}\}$ and identify the vertices\/ $x_{1}$ and\/ $x_{2}$ we obtain a graph\/ $X_{2}$ in which\/ $L_{1}$ and\/ $L_{2}$ are edge-equivalent but not vertex-equivalent. In other words, there exist two vertex-ends but only one edge-end, see Figure 4b.\/ $L_{1}$ and\/ $L_{2}$, and therefore their corresponding vertex-ends lie in every neighbourhood of\/ $x$ in the vertex-topology. Thus\/ $(VX_{2}\cup \Omega_{\mathrm{\scriptscriptstyle V}}X_{2},$\tauv$X_{2})$ is not a\/ $T_{1}$-space. In the edge-topology the edge-end of\/ $X_{2}$ and the vertex\/ $x$ cannot even be separated in the sense of the\/ $T_{0}$-axiom. In the graph\/ $X_{1}$ the same holds for the edge-end of $L_{1}$ and $L_{2}$ and the vertices\/ $x_{1}$ and\/ $x_{2}$, respectively.
\end{exmp}

\begin{center}
\setlength{\unitlength}{1mm}
\begin{picture}(110,65)
\put(0,0){\vector(0,1){55}}
\put(50,0){\vector(0,1){55}}
\put(70,0){\vector(0,1){55}}
\put(110,0){\vector(0,1){55}}
\put(70,0){\line(2,1){40}}
\put(70,10){\line(1,0){40}}
\put(70,20){\line(2,-1){40}}

\put(0,0){\line(2,1){20}}
\put(0,10){\line(1,0){50}}
\put(0,20){\line(2,-1){20}}

\put(50,0){\line(-2,1){20}}
\put(50,20){\line(-2,-1){20}}
\put(20,10){\line(1,0){10}}

\multiput(0,40)(8,-12){3}{\line(2,-3){4}}
\multiput(50,40)(-8,-12){3}{\line(-2,-3){4}}

\multiput(70,40)(8,-12){3}{\line(2,-3){4}}
\multiput(110,40)(-8,-12){3}{\line(-2,-3){4}}

\put(0,0){\circle*{1}}
\put(0,10){\circle*{1}}
\put(0,20){\circle*{1}}
\put(0,40){\circle*{1}}
\put(20,10){\circle*{1}}
\put(30,10){\circle*{1}}
\put(50,0){\circle*{1}}
\put(50,10){\circle*{1}}
\put(50,20){\circle*{1}}
\put(50,40){\circle*{1}}
\put(70,0){\circle*{1}}
\put(70,10){\circle*{1}}
\put(70,20){\circle*{1}}
\put(70,40){\circle*{1}}
\put(90,10){\circle*{1}}

\put(110,0){\circle*{1}}
\put(110,10){\circle*{1}}
\put(110,20){\circle*{1}}
\put(110,40){\circle*{1}}
\put(-2,60){$L_{1}$}
\put(48,60){$L_{2}$}
\put(68,60){$L_{1}$}
\put(108,60){$L_{2}$}
\put(88,5){$x$}
\put(20,5){$x_{1}$}
\put(30,5){$x_{2}$}
\put(23,30){$X_{1}$}
\put(88,30){$X_{2}$}

\end{picture}
\par

{\em \quad\, Figure 4a\hspace{46mm} Figure 4b}
\end{center}
\vspace{0.4mm}
As shown in Example 3 there may exist vertices with infinite degree, whose neighbourhoods in the
edge-topology all contain a fixed edge-end. To avoid such complications and to obtain better properties of separation the edge equivalence can be extended to the set of rays and vertices with infinite degree. In other words, vertices with infinite degree are considered as degenerated rays. They also can define new ends, so-called \emph{improper vertices}. We then obtain \emph{edge topologies \tauee\  of the second type} which are compact, totally disconnected and normal, either by defining $\{x\}$ as open sets a priori for all vertices $x$, see \cite{cartwright93martin} or by considering vertices $x$ with infinite degree as degenerated rays only and not as vertices (cf. \cite{kroen98topologische}). In the first case the edge-topology is a compactification of the discrete topology on $VX$. The compactness of these modified topologies can be deduced easily from Theorem 2.

\begin{thm}
Let\/ $c_{\mathrm{\scriptscriptstyle E}}=e_{\mathrm{\scriptscriptstyle E}}\cup\varepsilon_{\mathrm{\scriptscriptstyle E}}$ and\/ $c_{\mathrm{\scriptscriptstyle V}}=e_{\mathrm{\scriptscriptstyle V}}\cup\varepsilon_{\mathrm{\scriptscriptstyle V}}$ be such that\/ $e_{\mathrm{\scriptscriptstyle E}} \subset VX$, $e_{\mathrm{\scriptscriptstyle V}} \subset VX$, $\varepsilon_{\mathrm{\scriptscriptstyle E}}\subset \Omega_{\mathrm{\scriptscriptstyle E}}X$ and $\varepsilon_{\mathrm{\scriptscriptstyle V}}\subset \Omega_{\mathrm{\scriptscriptstyle V}}X$. $c_{\mathrm{\scriptscriptstyle E}}$ is open and closed in \taue\ if and only if\/ $e$ is an edge-cut and\/ $\varepsilon_{\mathrm{\scriptscriptstyle E}}=\Omega_{\mathrm{\scriptscriptstyle E}}e$. The set\/ $c_{\mathrm{\scriptscriptstyle V}}$ is open and closed in \tauv\ if and only if\/ $e$ is an edge-cut and $\varepsilon_{\mathrm{\scriptscriptstyle V}}=\Omega_{\mathrm{\scriptscriptstyle V}}e$.
\end{thm}

\begin{proof}
Let $c_{\mathrm{\scriptscriptstyle E}}$ and $c_{\mathrm{\scriptscriptstyle V}}$ be open and closed in their corresponding topologies. They can be represented as a finite union of elements of the base, because they are open and compact. Hence $e_{\mathrm{\scriptscriptstyle E}}$ is an edge-cut and $e_{\mathrm{\scriptscriptstyle V}}$ is a vertex-cut. The same argument holds for $(VX\cup\Omega_{\mathrm{\scriptscriptstyle V}}X)\backslash c_{\mathrm{\scriptscriptstyle V}}$ and therefore $e_{\mathrm{\scriptscriptstyle V}}^{*}$ is a vertex-cut, too. Thus $e_{\mathrm{\scriptscriptstyle V}}$ is an edge-cut.\par
An end in $\varepsilon_{\mathrm{\scriptscriptstyle E}}$ or $\varepsilon_{\mathrm{\scriptscriptstyle V}}$ must lie in $e_{\mathrm{\scriptscriptstyle E}}$ or $e_{\mathrm{\scriptscriptstyle V}}$, respectively, because $c_{\mathrm{\scriptscriptstyle E}}$ and $c_{\mathrm{\scriptscriptstyle V}}$ are open. On the other hand every end lying in $e_{\mathrm{\scriptscriptstyle E}}$ or $e_{\mathrm{\scriptscriptstyle V}}$ must be an element of $\varepsilon_{\mathrm{\scriptscriptstyle E}}$ or $\varepsilon_{\mathrm{\scriptscriptstyle V}}$, respectively. Otherwise the complements of $c_{\mathrm{\scriptscriptstyle E}}$ and $c_{\mathrm{\scriptscriptstyle V}}$ would not be open.
\end{proof}

\section{Metric ends}

Sometimes theorems for locally finite graphs can be generalized for non-locally-finite graphs or their proofs
can be simplified by replacing arguments that use the finiteness of sets of vertices by arguments of finite
diameters (cf. \cite{kroen98topologische}). At the other hand there even exist rather simple structures whose ramifications can neither be captured by the vertex- nor by the
edge-topology (cf. Example 4 and Example 5). This motivates a new definition of ends in non-locally-finite graphs, which is only based on the natural graph metric.\par

A \emph{metric cut} is a set of vertices with a vertex-boundary of finite diameter. A ray whose infinite
subsequences have all infinite diameters is called \emph{metric ray}. Two metric rays are \emph{metrically equivalent}, if they cannot be separated by metric cuts. Metrical equivalence is an equivalence relation on the set of metric rays of a graph. We call its
equivalence classes \emph{proper metric ends}. A proper metric  end \emph{lies in a set of vertices} $e$ if
all of its rays lie in $e$. 

We denote the set of proper metric ends that lie in $e$ by
$\Omega_{\mathrm{\scriptscriptstyle P} } e$ and write $\Omega_{\mathrm{\scriptscriptstyle P}} X$ instead
of $\Omega_{\mathrm{\scriptscriptstyle P}} VX$.

\begin{lem}
For any graph\/ $X$, the set\/ $B_{\mathrm{\scriptscriptstyle P}}X:=\{e\cup
\Omega_{\mathrm{\scriptscriptstyle P}}e \mid e \subset VX \mathit{\ and\ } \diam\theta e <
\infty \} $ is closed under finite intersection.
\end{lem}
The proof can be copied word by word from Lemma 2.
\par
We now can define the \emph{proper metric end topology} \taup$X$ as the topology on $VX\cup
\Omega_{\mathrm{\scriptscriptstyle P}}X$ which is generated by $B_{\mathrm{\scriptscriptstyle 
P}}X$.

\begin{thm}
The proper metric end topology is a totally disconnected Tychonoff topology. In the case that\/ $X$ is locally countably infinite it is Lindel\"of and normal. If there exists a vertex with infinite degree it is not paracompact.
\end{thm}
\begin{proof}
For any closed set $F$ not containing some point $x$ we can find a neighbourhood $A$ of $x$, which is element of the base and contained in the complement of $F$. As the base elements in $B_{\mathrm{\scriptscriptstyle P}}X$ are open and closed the indicator function on $A$ is continuous. Thus the proper metric end topology is $T_{3 \frac{1}{2}}$ or completely regular. Zero-dimensionality is also an immediate consequence of the fact that elements of the base are open and closed. It is easy to see that this topology is Hausdorff which now implies that it is Tychnoff and totally disconnected.\par
Let $X$ be a locally countable graph. For some given vertex $x$ let ${\mathcal V}$ denote the set of all base elements $e\cup\Omega_{\mathrm{\scriptscriptstyle P}}e$ such that $e$ is a connected component of $K(x,r)^{*}$ for some natural $r$. ${\mathcal V}$ is countable. Any element of an open cover ${\mathcal U}$ containing a proper metric end must contain an element of ${\mathcal V}$. Thus we can find a countable subcover of ${\mathcal U}$ of the set of proper metric ends. The set of vertices is countable anyway and hence the proper metric end topology is Lindel\"of. Every regular Lindel\"of space is normal.\par

For any vertex $x$ with infinite degree
$$\{\{x,y\}\mid x\sim y\}\cup\{K(x,1)^{*}\cup\, \Omega_{\mathrm{\scriptscriptstyle P}}X\}$$
is an open cover which is not locally finite in $x$ and therefore the proper metric end topology is not paracompact.
\end{proof}

We do not know whether the proper metric end topology is normal in the general case.

\begin{exmp}
Let\/ $X$ be the graph which arises from a tree\/ $T$ with only vertices of infinite degree to which we add edges
connecting all pairs of vertices with distance two. The proper metric topologies on\/ $X$ and\/ $T$ are
homeomorphic, in contrast to the vertex- and the edge-topology on\/ $X$, which are indiscrete and therefore
describe no structure at all.
\end{exmp}
To motivate the definitions of the following section we remark that every sequence of vertices
that has neither an accumulation point in\/ $VX$ (equivalently: has infinitely many identical elements) nor an
accumulation point in\/ $\Omega_{\mathrm{\scriptscriptstyle P}}X$ can be divided up into two types of
subsequences:

\begin{enumerate}
\item Bounded sequences.
\item Sequences with no bounded subsequences for which there exists a metric cut, whose complement consists
only of components that contain at most finitely many elements of the sequence. 
\end{enumerate}

\section{A compactification of the proper metric end topology}

We will now add two additional points to the set of metric ends and modify the proper metric topology
correspondingly.\par
A \emph{star-cut} is a metric cut $e$ such that every union $r$ of all but finitely many components of its
complement has infinite diameter. For a star-cut $e$ we call such a union $r$ \emph{star-boundary of}
$e$. 
A set $s$ of vertices is called a \emph{global star-set} if to every star-cut $e$ there exists a
star-boundary $r$, such that $r\backslash s$ has a finite diameter.\par
An infinite set of vertices $p$ is called \emph{locally complete}, if it contains all but finitely many
elements of every ball.\par
If there exists a star-cut in $VX$ we add an element $\sigma_{X}$ called the \emph{star-end} to the set of
proper metric ends. Furthermore we add an element $\lambda$ to the set of ends which we call the \emph{local end}. The set $\Omega_{\mathrm{\scriptscriptstyle M}}X$ so obtained is called \emph{set of metric ends of} $X$, $\lambda_{X}$ 
and $\sigma_{X}$ are called \emph{improper metric ends}. We say
that the local end \emph{lies in a set of vertices} $e$, if it is a locally complete set. If $e$ is a global star-set we say that the star-end \emph{lies in} $e$. Let $\Omega_{\mathrm{\scriptscriptstyle M}}e$ denote the set of all metric ends, that lie in $e$.\par
Although the local end will be useless in locally finite graphs (equivalently: it is an isolated point, $\{\lambda_{X}\}$ is open), we add it to the set of ends in any way. This will help use to find
homeomorphisms between end-spaces in locally and non-locally-finite graphs.
\begin{lem}
The set 
$$B_{\mathrm{\scriptscriptstyle M}}X:=\{(e\cap s\cap p)\cup
\Omega_{\mathrm{\scriptscriptstyle M}}(e\cap s\cap p) \mid e\text{ is a metric cut, }$$
$$s\text{ is a global star-set and }p\text{ is a locally complete set,}\}$$
is closed under finite intersection.
\end{lem}

We call the topology \taum$X$ on $VX\cup\Omega_{\mathrm{\scriptscriptstyle M}}X$, which is generated by
$B_{\mathrm{\scriptscriptstyle M}}X$, \emph{metric end topology}.

\begin{lem}
$(VX\cup\Omega_{\mathrm{\scriptscriptstyle M}}X,$\taum$X)$ is a Lindel\"of space.
\end{lem}

\begin{proof}
Let $\mathcal{U}$ be an open cover of base elements. It must contain sets
$p\cup\Omega_{\mathrm{\scriptscriptstyle M}}p$ and $s\cup\Omega_{\mathrm{\scriptscriptstyle
M}}s$ such that $p$ is a locally complete set and $s$ is a global star-set. After removing these
base elements from $VX\cup\Omega_{\mathrm{\scriptscriptstyle M}}X$ we obtain a set
$A\cup\Omega A$ such that $A\subset p^{*}\cap s^{*}$ and $\Omega A\subset 
\Omega_{\mathrm{\scriptscriptstyle
P}}X\backslash\Omega_{\mathrm{\scriptscriptstyle P}}s$. $A$ is a set of
vertices, that contains only finitely many elements of any ball in $X$. Thus $A$ is countable and we can find a countable subcover $\mathcal{V}$ of $VX$ in $\mathcal{U}$. We now copy
the ideas in the second part of the proof of Lemma 6 to construct a countable subset of
$\mathcal{U}$ which covers $\Omega A$.\par
Let $\dm_{\mathcal{\, U}}(\omega,z)$ again denote the distance of the end $\omega$ to the vertex $z$
with respect to the cover $\mathcal{U}$. As $A$ contains no
star-boundary, it has the property that for every metric cut $e$ we can find a union $u$ of finitely many
components of $e^{*}$, such that $A\backslash u$ has a finite diameter. 
Thus there are only finitely many connected components in $(s\cup
K(z,n))^{*}$ for any natural number $n$, we can find a finite subcover $\mathcal{V}_{n}$ of $\mathcal{U}$
covering all ends with $\dm_{\mathcal{\, U}}(\omega,z)\le n$, that lie in $\Omega A$. Thus
$$\{p\cup\Omega_{\mathrm{\scriptscriptstyle M}}p\}\cup\{s\cup\Omega_{\mathrm{\scriptscriptstyle
M}}s\}\cup\mathcal{V}\cup\bigcup_{n=1}^{\infty}\mathcal{V}_{n}$$
is a countably finite subcover of $\mathcal{U}$ covering $VX\cup \Omega_{\mathrm{\scriptscriptstyle
M}}X$.
\end{proof}

\begin{thm}
$(VX\cup\Omega_{\mathrm{\scriptscriptstyle M}}X,$\taum$X)$ is a compact topological space. All limits of
convergent sequences are unique. For a convergent sequence\/ $\xi$ with\/ $\alpha =\lim \xi$ exactly one of the
following cases must hold:

\begin{enumerate}
\item The limit\/ $\alpha$ is a vertex. All elements of the sequence are equal to\/ $\alpha$ from an index on.
\item From some index onwards the elements of the sequence\/ $\xi$ are the local end or they are vertices lying in some ball such that at the most finitely many vertices are identical. In this case\/ $\xi$ converges to the local end\/ $\lambda_{X}$.
\item In the complement of any ball exactly one component contains all but finitely many elements
of the sequence $\xi$. There exists exactly one proper metric end lying in all these components, which is the
limit\/ $\alpha$.
\item There exists a vertex\/ $x$ and a radius\/ $r$ such that for any natural\/ $n$ the ball\/ $K(z,r+n)$ is a
star-cut, whose star-boundaries contain all but finitely many of those elements of the sequence, that do not
equal the star-end, but contain only finitely many of them in each of their connected components. The
sequence converges to the star-end\/ $\sigma_{X}$.
\end{enumerate}
\end{thm}

\begin{proof}
To show that every sequence in $VX\cup\Omega_{\mathrm{\scriptscriptstyle M}}X$ has an 
accumulation point we can assume that it contains an infinite partial sequence $\xi$ consisting of pairwise
different elements. Otherwise the existence of an accumulation point is immediate. If $\xi$ contains a
bounded subsequence of vertices, the local end is an accumulation point. In the other case we choose a vertex
$x$ and distinguish between two cases.
\begin{enumerate}
\item For every natural $n$ there exists a sequence  $(e_{n})_{n \in \mathbb{N}}$ of components of
$K(x,n)^{*}$ such that $e_{n}$ contains infinitely many elements of the sequence $\xi$ and is a superset
of $e_{n+1}$. Following the idea of Lemma 3 we now can construct a metric ray, that lies in all cuts $e_{n}$.
Its end is an accumulation point of $\xi$.
\item There exists a ball $K(x,n)$ such that only finitely many elements of $\xi$ lie in every component of
$K(x,n)^{*}$. Let $s\cup \Omega_{\mathrm{\scriptscriptstyle M}}s$ be a base element that contains the
star-end $\sigma_{X}$, which means that $s$ must be a global star-set. By the definition of a global star-set
there exists a star-annulus $r$ of $K(x,n)$ such that $r\backslash s$ has a finite diameter. The set $r\cup
\Omega_{\mathrm{\scriptscriptstyle M}}r$ must contain all but finitely many elements of $\xi$, but as $\xi$ 
contains no bounded subsequence, only finitely many elements of the sequence $\xi$ lie in $r\backslash s$.
Now $s\cup \Omega_{\mathrm{\scriptscriptstyle M}}s$ must contain all but finitely many elements of $\xi$,
which means that $\sigma_{X}$ in an accumulation point of $\xi$.
\end{enumerate}
We have shown that $(VX\cup\Omega_{\mathrm{\scriptscriptstyle M}}X,$\taum$X)$ is sequentially compact. By
Lemma 10 it is also a Lindel\"of space and hence it is compact. The other statements of Theorem 5 now follow
easily by the above considerations.
\end{proof}

\begin{exmp}
Let\/ $K_{\mathbb{N}}$ denote the complete graph with vertex-set\/ $\mathbb{N}$, see Figure 5a. We take two graphs\/ $K_{\mathbb{N}}^{(1)}$ and\/ $K_{\mathbb{N}}^{(2)}$ that are isomorphic to\/ $K_{\mathbb{N}}$ and connect each vertex with adjacent vertices\/ $x_{1}$ and\/ $x_{2}$, respectively (Figure 5b). As in Example 3 we take away the edge\/ $(x_{1},x_{2})$ from the graph in Figure 5b, identify\/ $x_{1}$ and $x_{2}$ to a vertex $x$. Then we connect each vertex in\/ $K_{\mathbb{N}}^{(1)}$ and\/ $K_{\mathbb{N}}^{(2)}$ with $x$ as shown in Figure 5c. In Figure 5d we take a sequence\/ $(K_{\mathbb{N}}^{(n)})_{n\in \mathbb{Z}}$ of copies of the graph\/ $K_{\mathbb{N}}$ and connect each vertex in\/ $K_{\mathbb{N}}^{(n)}$ with each of the vertices in\/ $K_{\mathbb{N}}^{(n+1)}$.\par
Every sequence of pairwise distinct vertices in the metric end topology of the graphs in Figure 5a, 5b and 5c converges to the pointend. For the edge-end and vertex-end compactification in Figure 5b and 5c we have the same situation as in Figure 4. In Figure 5d they contain only one end whereas the metric topology contains two proper metric ends, compare with Example 2.
\end{exmp}

\begin{center}

\setlength{\unitlength}{1mm}
\begin{picture}(115,45)
\put(10,35){\circle{20}}
\put(33,35){\circle{20}}
\put(42.98,35){\circle*{1}}
\put(42.98,35){\line(-1,1){4.97}}
\put(42.98,35){\line(-1,-1){4.97}}
\put(47.98,35){\line(1,1){4.97}}
\put(47.98,35){\line(1,-1){4.97}}
\put(42.98,35){\line(1,0){5}}
\put(47.98,35){\circle*{1}}
\put(57.88,35){\circle{20}}
\put(80,35){\circle{20}}
\put(89.98,35){\line(-1,1){4.97}}
\put(89.98,35){\line(-1,-1){4.97}}
\put(89.98,35){\line(1,1){4.97}}
\put(89.98,35){\line(1,-1){4.97}}
\put(89.98,35){\circle*{1}}
\put(99.88,35){\circle{20}}
\multiput(25,0)(15,0){5}{\line(0,1){15}}
\put(25,0){\vector(-1,0){15}}
\put(25,0){\vector(1,0){75}}
\put(25,15){\vector(-1,0){15}}
\put(25,15){\vector(1,0){75}}
\put(3.5,22){{\em Fig.\ 5a}}
\put(39,22){{\em Fig.\ 5b}}
\put(83,22){{\em Fig.\ 5c}}
\put(7.5,34){$K_{\mathbb{N}}$}
\put(29.5,34){$K_{\mathbb{N}}^{(1)}$}
\put(54.38,34){$K_{\mathbb{N}}^{(2)}$}
\put(76.5,34){$K_{\mathbb{N}}^{(1)}$}
\put(96.38,34){$K_{\mathbb{N}}^{(2)}$}
\multiput(29.5,6.5)(15,0){4}{$K_{\mathbb{N}}$}
\end{picture}
\par
\vskip2mm
{\em Figure 5d}
\end{center}

\section{Metric ends and quasi-isometries}

\begin{defn}
Two graphs $X$ and $Y$ are called \emph{quasi-isometric with respect to the functions} $\phi: VX\to VY$ \emph{and} $\psi: VY\to VX$ if there exist constants $a$, $b$, $c$ and $d$ such that for all vertices $x$,
$x_{1}$ and $x_{2}$ in $VX$ and vertices $y$, $y_{1}$ and $y_{2}$ in $VY$, the following conditions hold

\vskip6pt
\begin{tabular}{lll}
(Q1) & $\dm_{Y}(\phi(x_{1}),\phi(x_{2}))\le a\cdot\dm_{X}(x_{1}, x_{2})$\qquad & (boundedness of $\phi$)\\
(Q2) & $\dm_{X}(\psi(y_{1}),\psi(y_{2}))\le b\cdot\dm_{Y}(y_{1}, y_{2})$ & (boundedness of $\psi$)\\
(Q3) & $\dm_{X}(\psi\phi(x),x)\le c$ & (quasiinjectivity of $\phi$)\\
(Q4) & $\dm_{Y}(\phi\psi(y),y)\le d$ & (quasisurjectivity of $\phi$)
\end{tabular}
\vskip6pt
We call $\phi$ and $\psi$ \emph{quasi-inverse} to each other.
\end{defn}
For general metric spaces the definition of quasi-isometries allows further additive constants in the Axioms (Q1) and (Q2). In case the positive values of the metric are greater than some positive real number these additive constants are useless.\par
Quasi-isometries may change structures as long as the differences can be bounded uniformly. In other words we could say that they preserve the global structure of graphs when we consider graphs as discrete metric spaces only.\par
Without proof we remark that quasi-isometry is an equivalence relation on the family of all graphs.

\begin{exmp}
The graphs in Example 3, the ladder graph (Example 2) and the graph in Example 5, Figure 4d, as well as the tree\/ $T$ and the graph\/ $X$ in Example 4 are quasi-isometric.\par
For all constants\/ $r$ the Cayley graphs\/ $X$ of the free group with respect to the generating system\/ $A^{r}$ in Section 8 are quasi-isometric to each other.\par
Replacing every edge in a graph $X$ by a path with a length smaller or equal some constant we obtain a graph which is quasi-isometric to\/ $X$.
\end{exmp}

\begin{lem}
Let\/ $\phi: VX\to VY$ be a quasi-isometry and\/ $A$ a subset of\/ $VX$, then
$$\diam_{X}A<\infty \Leftrightarrow \diam_{Y}\phi(A)<\infty.$$
\end{lem}

\begin{proof}
$\diam_{X}A<\infty$ implies
$$\diam_{Y}\phi(A)=\{\dm_{Y}(\phi(x_{1}),\phi(x_{2}))\mid x_{1}, x_{2}\in A\}\le$$
$$\sup\{ a\cdot\dm_{X}(x_{1}, x_{2})\mid x_{1}, x_{2}\in A\}=a\cdot\diam_{X}A<\infty.$$
Let $\psi$ be a quasi-inverse to $\phi$. If $\diam_{X}A$ is infinite, then, by (Q3), this also holds for
$\diam_{X}\psi\phi (A)$. Now $$\dm_{X}(\psi\phi(x_{1}),\psi\phi(x_{2}))\le b\cdot\dm_{Y}(\phi(x_{1}),
\phi(x_{2}))$$ implies that $$\{\dm_{Y}(\phi(x_{1}),\phi(x_{2}))\mid x_{1},x_{2}\in A\}$$ has no upper
bound.
\end{proof}

\begin{cor}
The pre-image of a metric cut or a global star-set set under a quasi-isometry is a metric cut or a global star-set, respectively. 
\end{cor}

We now want to extend the concept of quasi-isometry to the set of proper metric ends.

\begin{thm}
To every quasi-isometry\/ $\phi: VX\to VY$ there exists a unique extension
$$\Phi:VX\cup\Omega_{\mathrm{\scriptscriptstyle P}}X\to VY\cup\Omega_{\mathrm{\scriptscriptstyle P}}Y,$$
such that
\begin{enumerate}
\item $\Phi|_{VX}=\phi$
\item $\Phi$ is continuous and
\item $\Phi|_{\Omega_{\mathrm{\scriptscriptstyle P}}X}$ is a homeomorphism of
$\Omega_{\mathrm{\scriptscriptstyle P}}X$ and\/ $\Omega_{\mathrm{\scriptscriptstyle P}}Y$
with respect to the corresponding relative topologies.
\end{enumerate}
\end{thm}

\begin{proof}
By connecting the $\phi$-images of adjacent vertices of a metric ray $L$ in $X$ with geodesic paths with
lengths, that are smaller or equal $a$, we obtain a path $P$ in $Y$. Its diameter is infinite by Lemma 11.
If it had infinite subset $M$ with finite diameter, we also could find infinitely many elements in $M$
that lie in $\phi(L)$, contradicting Lemma 11 and the assumption that $L$ is a metric ray. As a graph is locally finite if and only if every bounded set of vertices is finite, $P$ is a locally finite subgraph of $Y$. Lemma 1 implies that it must contain a ray which we denote with
$\tilde{\phi}(L)$. By the above consideration it is also a metric ray in $Y$. Thus $\tilde{\phi}$ maps metric rays in $X$ onto metric rays in $Y$.\par 
Let $L_{1}$ and $L_{2}$ be metric equivalent rays in $X$. If $\tilde{\phi}(L_{1})$ and
$\tilde{\phi}(L_{2})$ were not metrically equivalent, we could find disjoint metric cuts $f_{1}$ and
$f_{2}$, such that $\tilde{\phi}(L_{1})$ lies in $f_{1}$ and $\tilde{\phi}(L_{2})$ lies in $f_{2}$. Again
by Lemma 10,
$\phi^{-1}(f_{1})\cap\phi^{-1}(f_{2})$ must have a finite diameter and hence the rays $L_{1}$ and $L_{2}$
cannot lie in both cuts $f_{1}$ and $f_{2}$ in contradiction to the assumption that they are metrically
equivalent. Thus we can say that metric equivalence is an invariance under $\tilde{\phi}$ on sets of
metric rays.\par
Now for every end $\omega$ in $\Omega_{\mathrm{\scriptscriptstyle P}}X$ we define $\Phi(\omega)$ as the
unique  end in $\Omega_{\mathrm{\scriptscriptstyle P}}Y$ which contains the $\tilde{\phi}$\/-images of the
elements of $\omega$ and set $\Phi(x)=\phi(x)$ for every vertex $x$.\par

By the above invariance of the metric equivalence and the construction of $\tilde{\phi}$ we obtain
$$\Phi (\Omega_{\mathrm{\scriptscriptstyle P}}e)\subset \Omega_{\mathrm{\scriptscriptstyle P}} \Phi(e)$$
for every metric cut $e$ in $X$, and therefore by Corollary 1
$$\Phi (\Omega_{\mathrm{\scriptscriptstyle P}}\Phi^{-1}(f))\subset \Omega_{\mathrm{\scriptscriptstyle P}}
\Phi(\Phi^{-1}(f))=\Omega_{\mathrm{\scriptscriptstyle P}}f$$
and
$$\Omega_{\mathrm{\scriptscriptstyle P}}\Phi^{-1}(f)\subset \Phi^{-1}(\Omega_{\mathrm{\scriptscriptstyle P}}
f)$$
for every metric cut $f$ in $Y$.\par
Now let $\omega$ be a proper metric end in $\Phi^{-1}(\Omega_{\mathrm{\scriptscriptstyle P}}f)$. Then of
course $\Phi(\omega)\in \Omega_{\mathrm{\scriptscriptstyle P}} f$. If $\omega$ did not lie in
$\Phi^{-1}(f)$ then there would exist infinitely many vertices $x$ in a ray $L$ of $\omega$ which are not
elements of $\Phi^{-1}(f)$. Their $\phi$-images would not be elements of $f$, $\tilde{\phi}(L)$ would not lie in $f$ and $\Phi(\omega)$ would not element of $\Omega_{\mathrm{\scriptscriptstyle P}} f$. Hence
$$\Phi^{-1}(\Omega_{\mathrm{\scriptscriptstyle P}}f)\subset \Omega_{\mathrm{\scriptscriptstyle
P}}\Phi^{-1}(f)$$
and

\begin{equation}
  \label{eq:hbeta1}
\Phi^{-1}(\Omega_{\mathrm{\scriptscriptstyle P}}f)=\Omega_{\mathrm{\scriptscriptstyle P}}\Phi^{-1}(f).
\end{equation}

For every base element $f\cup \Omega_{\mathrm{\scriptscriptstyle P}}f$ in $B_{\mathrm{\scriptscriptstyle
P}}Y$ we obtain
$$\Phi^{-1}(f\cup \Omega_{\mathrm{\scriptscriptstyle P}}f)=\Phi^{-1}(f)\cup
\Omega_{\mathrm{\scriptscriptstyle P}} \Phi^{-1}(f)$$
and $\Phi$ is continuous.\par
To prove the third property of $\Phi$, it now suffices to show that the restriction
$\Phi|_{\Omega_{\mathrm{\scriptscriptstyle P}}X}$ is a bijection of the sets of ends in $X$ and $Y$.
To a given end $\varepsilon$ in $\Omega_{\mathrm{\scriptscriptstyle P}}Y$ we can find a decreasing sequence 
$(f_{n})_{n\in \mathbb{N}}$ of connected metric cuts in $Y$ such that 
$\bigcap_{n\in \mathbb{N}}f_{n}=\emptyset$ and $\bigcap_{n\in \mathbb{N}} \Omega_{\mathrm{\scriptscriptstyle
P}} f_{n}=\varepsilon$. Now we choose an increasing sequence $(K(z,r_{n}))_{n\in \mathbb{N}}$ of concentric
balls in $X$ such that $\theta \Phi^{-1}(f_{n})$ is contained in $K(z,r_{n})$. To every natural $n$
there exists a connected component $e_{n}$ of $K(z,r_{n})^{*}$ which contains all but finitely many
vertices of the $\Phi$-preimage of a ray $L$ of $\varepsilon$. Copying the idea of Lemma 3 we can find a
metric ray $R$ that lies in all cuts $e_{n}$. The end $\omega$ of $R$ is mapped onto $\varepsilon$ under
$\Phi$. Thus $\Phi$ is surjective.\par
Let $\varepsilon_{1}$ and $\varepsilon_{2}$ be two different ends of the graph $Y$ and $f_{1}$ and $f_{2}$
two disjoint metric cuts such that $\varepsilon_{1}$ lies in $f_{1}$ and $\varepsilon_{2}$ lies in $f_{2}$.
Now $\Phi^{-1}(\Omega_{\mathrm{\scriptscriptstyle P}} f_{1})$ and
$\Phi^{-1}(\Omega_{\mathrm{\scriptscriptstyle P}} f_{2})$ must be disjoint too, and therefore $\Phi$ is
injective.\par
To complete the proof of the theorem, it remains to show that $\Phi$ is the unique extension of $\phi$ with
the requested properties. We assume that there exists another extension $\hat{\Phi}$ which has these three
properties but does not equal $\Phi$. Let $\omega$ be an end in $X$ and $\Phi(\omega)=\varepsilon_{1}$ and
$\hat{\Phi}(\omega)=\varepsilon_{2}$ for two different ends $\varepsilon_{1}$ and $\varepsilon_{2}$ in $Y$.
Now we choose again two disjoint metric cuts $f_{1}$ and $f_{2}$ such that $\varepsilon_{1}$ lies in $f_{1}$
and $\varepsilon_{2}$ lies in $f_{2}$. As $\Phi^{-1}(f_{1})=\phi^{-1}(f_{1})$ and
$\hat{\Phi}^{-1}(f_{2})=\phi^{-1}(f_{2})$, the preimages $\Phi^{-1}(f_{1})$ and $\hat{\Phi}^{-1}(f_{2})$ are
disjoint metric cuts in $X$. We know that $\omega$ must 
lie in $\Phi^{-1}(f_{1})$ and therefore $\omega$ must lie in $\hat{\Phi}^{-1}(f_{2})^{*}$. This implies that
every open neighbourhood of $\omega$ has a nonempty intersection with $\hat{\Phi}^{-1}(f_{2})^{*}$. In other
words, no open neighbourhood of $\omega$ is completely contained in $\hat{\Phi}^{-1}(f_{2}\cup
\Omega_{\mathrm{\scriptscriptstyle P}} f_{2})$ and $\hat{\Phi}$ is not continuous.
\end{proof}

For a subset $A$ of $VX\cup \Omega_{\mathrm{\scriptscriptstyle P}}X$ and a
natural number $r$ we define
$$A+r:=\{x\in VX\mid\dm_{X}(x,A\backslash\Omega_{\mathrm{\scriptscriptstyle P}}A
)\le r\}\cup \Omega_{\mathrm{\scriptscriptstyle P}}A.$$
In the following sense we could call $\Phi$ a {\em quasi-open\/}
function.
\begin{thm}
For any base element\/ $e\cup \Omega_{\mathrm{\scriptscriptstyle P}} e$ in\/ $B_{\mathrm{\scriptscriptstyle
P}}X$, the set\/ $\Phi(e\cup \Omega_{\mathrm{\scriptscriptstyle P}} e)+d+1$ is open 
in \taup$Y$.
\end{thm}

Although the proof seems a little technical, its idea is simple. Quasi-isometry is a weakened form of
isomorphy. The quasi-surjectivity does not ensure that the image $\phi(VX)$ covers $VY$ 
completely but says that it does not have `holes' that are bigger than a circle with radius $d$.
\begin{proof}
Let $y$ be an element of $\theta(\phi(e)+d+1)$. As $\dm_{Y}(\phi\psi(y),y)\le d$ the vertex $\phi\psi(y)$
must not be contained in $\phi(e)$, so $\phi\psi(y)\in \phi(e^{*})$ and $\psi(y)\in e^{*}$. Hence
$$\psi(\theta(\phi(e)+d+1))\subset e^{*}.$$
Now let $y_{1}$ be an element of $\phi(e)+d+1$ which is adjacent to $y$ and choose an element $y_{2}$ of
$\phi(e)$ such that
$\dm_{Y}(y_{1},y_{2})\le d+1.$
For an $x$ in $e$ with $\phi(x)=y_{2}$ we now have
$$\dm_{Y}(y,\phi(x))\le d+2.$$ By the boundedness of $\psi$ we have 
$\dm_{X}(\psi(y),\psi\phi(x))\le b\cdot(d+2)$
and by the quasi-injectivity of $\phi$ we get 
$\dm_{X}(x,\psi\phi(x))\le c$ which implies
$$\dm_{X}(\psi(y),x)\le b\cdot (d+2)+c.$$ In other words
$$\max \{\dm_{X}(\psi(y),e)\mid y\in \theta(\phi(e)+d+1) \}\le  b\cdot (d+2)+c.$$

As $\psi(\theta(\phi(e)+d+1))$ is a subset of $e^{*}$ we now obtain
$$\psi(\theta(\phi(e)+d+1))\subset \theta e+b\cdot (d+2)+c-1.$$
Consequently $\psi(\theta(\phi(e)+d+1))$ has a finite diameter and by Lemma 11 this also holds for
$\theta(\phi(e)+d+1)$ and $\phi(e)+d+1$ is a metric cut in $Y$.\par
By the definition of $\Phi$
$$\Phi(e\cup \Omega_{\mathrm{\scriptscriptstyle P}} e)+d+1=\Phi(\Omega_{\mathrm{\scriptscriptstyle P}}
e)\cup (\phi(e)+d+1)$$
and to prove that this set is open in \taup$Y$ we show that $\Phi(\Omega_{\mathrm{\scriptscriptstyle P}}e)$
is contained 
in $\Omega_{\mathrm{\scriptscriptstyle P}}(\phi(e)+d+1)$. For every metric cut $f$ in $Y$ the equation (7.1)
implies 
$$\Omega_{\mathrm{\scriptscriptstyle P}} f=\Phi(\Omega_{\mathrm{\scriptscriptstyle P}}\Phi^{-1}(f)).$$
As $\Omega_{\mathrm{\scriptscriptstyle P}}e$ is a subset of $\Omega_{\mathrm{\scriptscriptstyle
P}}\Phi^{-1}(\phi(e)+d+1)$ and $\phi(e)+d+1$ is a metric cut we now obtain
$$\Phi(\Omega_{\mathrm{\scriptscriptstyle P}}e)\subset \Phi(\Omega_{\mathrm{\scriptscriptstyle
P}}\Phi^{-1}(\phi(e)+d+1))=\Omega_{\mathrm{\scriptscriptstyle P}}(\phi(e)+d+1)$$
by replacing $f$ with $\phi(e)+d+1$.

\end{proof}

For applications of quasi-isometries in the study of ends of graphs confer \cite{moeller92ends2} and \cite{kroen98topologische}.

\section{Bounded random walk on the free group}

We now want to give a further example for an application of the metric end compactification concerning the
random walk on the free group with countably infinitely many generators. Random walk on free groups were first introduced and studied by Kesten, \cite{kesten59symmetric}. Let $\Gamma$ be a free group with a symmetric and countable set of generators $A$ containing at least four elements. Every element $x$ of $\Gamma$ can be represented in a unique way by the shortest product of elements of $A$ that equals $x$. The length of this product is called \emph{length} of $x$, the length of the neutral element 
$o$ is set zero. Let $A^{r}$ denote the set of elements with positive lengths that are less or equal
some natural $r$. The Cayley graph $X$ of $\Gamma$ with respect to $A^{r}$ has vertex set
$VX=\Gamma$. 
Two vertices $x$ and $y$ are adjacent if and only if $x^{-1}y$ is an element of $A^{r}$. $X$ is isomorphic to a graph that arises from a regular tree $T$ with degree $|A|$ by connecting pairs of distinct vertices $x$ and $y$ with distance $\dm_{T}(x,y)\le r$. Compare with Example 4.\par
Let $\mu$ be a probability measure on $A^{r}$ whose support generates the whole group $\Gamma$. Now $\mu$ defines a random walk $(Z_{n})_{n\in \mathbb{N}}$ on
$\Gamma$ with respect to the probability measure $\p$ on $\Gamma^{\mathbb{N}}$ which is generated by the  cylindric sets in $(A^{r})^{\mathbb{N}}$ together with the corresponding powers of $\mu$.

\begin{thm}
The random walk\/ $(Z_{n})_{n\in \mathbb{N}}$ on the countably infinite free group\/ $\Gamma$ converges almost surely to some proper metric end in the Cayley graph\/ $X$.
\end{thm}

\begin{proof}
As the support of $\mu$ generates $\Gamma$ we can choose four distinct elements $a_{1}$, $a_{2}$,
$a_{1}^{-1}$ and $a_{2}^{-1}$ in $A^{r}$ such that $\mu(a_{1})$ and $\mu(a_{2})$ are positive. We define a random walk
$(\tilde{Z}_{n})_{n\in \mathbb{N}}$ on the free group 
$\tilde{\Gamma}$ with generating system $\tilde{A}=\{a_{1},a_{2},a_{1}^{-1},a_{2}^{-1}\}$ by choosing $\mu(a_{1})$,
$\mu(a_{2})$, $\mu(a_{1}^{-1})$ and $\mu(a_{2}^{-1})$ as the probabilities for the right multiplication with
the generating elements, respectively. The probability of not making a move is set $1-(\mu(a_{1})+\mu(a_{2})+\mu(a_{1}^{-1})+\mu(a_{2}^{-1}))$. We know that
$(\tilde{Z}_{n})_{n\in 
\mathbb{N}}$ is transient (cf. for example \cite{woess00random}). The corresponding probability measure on the set of
trajectories is denoted by $\tilde{\p}$. As the corresponding Cayley graph
$\tilde{X}$ is locally finite, the random walk $(\tilde{Z}_{n})_{n\in\mathbb{N}}$ does not enter any ball
from an index on with probability one. Thus
$$\tilde{\p}[\lim_{n\to \infty}\dm_{\tilde{X}}(o,\tilde{Z}_{n})=\infty]=1.$$
Now
$\dm_{\tilde{X}}(x,y)\le r\cdot \dm_{X}(x,y)$
for all elements $x$ and $y$ of $\Gamma$ and therefore
$$\p[\lim_{n\to \infty}\dm_{X}(o,Z_{n})=\infty]=1.$$
This is equivalent to the almost sure convergence of $Z_{n}$ to a proper metric end in the
metric end topology of $X$.
\end{proof}

\begin{appendix}
\vskip12pt
Section 2, 3 and 4 are part of the author's masters thesis \cite{kroen98topologische} at the University of Salzburg under the
supervision of 
Prof.\ W.\ Woess and the author wants to thank him for many useful suggestions. The main part of this thesis was written during a stay at Milan supported by the Italian
Ministry of Foreign Affairs.
\end{appendix}

\end{document}